\documentclass[10pt,a4paper,reqno]{amsart} 
\usepackage{epsfig,color}
\usepackage{graphicx}
\usepackage[english]{babel} 
\usepackage{amsmath}
\usepackage{amsthm,amssymb,latexsym} 
\usepackage{t1enc}
\usepackage{fancyheadings} 
\usepackage{epic} 
\usepackage{eepic}

\newtheorem{prop}{Proposition} 
\newtheorem{lemma}[prop]{Lemma}
\newtheorem{cor}[prop]{Corollary} 
\newtheorem{theorem}[prop]{Theorem}
\newtheorem{thm}[prop]{Theorem}

\theoremstyle{definition} 
\newtheorem{defin}[prop]{Definition}

\newcommand{\B}{\ensuremath{\mathcal{B}}} 
\newcommand{\X}{\ensuremath{\mathcal{X}}} 
\newcommand{\Y}{\ensuremath{\mathcal{Y}}} 

\newcommand{\A}{\ensuremath{\mathcal A}} 
\newcommand{\I}{\ensuremath{\mathcal I}} 
\newcommand{\ov}{\overline}
\newcommand{\un}{\underline}
\def\newop#1{\expandafter\def\csname #1\endcsname{\mathop{\rm#1}\nolimits}}
\newop{inv}

\newcommand{\aseq}{\A-sequence}
\newcommand{\bseq}{\B-sequence}

\newcommand{\bwx}{\phi}
\newcommand{\BWX}{\psi}
\newcommand{\Sn}{\ensuremath{\mathcal S}} 
\newcommand{\p}{permutation}
\newcommand{\ps}{permutations}
\newcommand{\gf}{generating function}
\newcommand{\gfs}{generating functions}
\newcommand{\lseq}{$\ell$-sequence}

\def\emm#1,{{\em #1}}
\def\lb{\linebreak}

\newcommand{\Ref}[1]{(\ref{#1})}

\catcode`\@=11
\def\section{\@startsection{section}{1}%
 \z@{.7\linespacing\@plus\linespacing}{.5\linespacing}%
 {\normalfont\bfseries\scshape\centering}}
\def\subsection{\@startsection{subsection}{2}%
  \z@{.5\linespacing\@plus\linespacing}{.5\linespacing}%
  {\normalfont\bfseries\scshape}}

\def\subsubsection{\@startsection{subsubsection}{3}%
  \z@{.5\linespacing\@plus.7\linespacing}{-.5em}%
  {\normalfont\itshape}}
\catcode`\@=12
\title[Decreasing subsequences and Wilf equivalence]{Decreasing
subsequences in permutations and Wilf equivalence for involutions}

\author{Mireille Bousquet-M\'elou \and Einar Steingr\'{\i}msson}
\address{CNRS, LaBRI, Universit\'e Bordeaux 1,
351 cours de la Lib\'eration,
  33405 Talence Cedex, France
\and   Matematik, Chalmers tekniska h\"ogskola och
  G\"oteborgs universitet  
S-412~96 G\"oteborg, Sweden}
\email{mireille.bousquet@labri.fr, einar@math.chalmers.se}
\thanks{Both authors were partially supported by the European Commission's IHRP
  Programme, grant HPRN-CT-2001-00272, ``Algebraic Combinatorics in
  Europe''}

\date{\today}

\begin{document}

\begin{abstract}
  In a recent paper, Backelin, West and Xin describe a map $\phi ^*$
  that recursively replaces all occurrences of the pattern $k\cdots
  21$ in a permutation $\sigma$ by occurrences of the pattern
  $(k-1)\cdots 21 k$. The resulting permutation $\phi^*(\sigma)$
  contains no decreasing subsequence of length $k$. We prove that,
  rather unexpectedly, the map $\phi ^*$ commutes with taking the
  inverse of a permutation.
  
  In the BWX paper, the definition of $\phi^*$ is actually extended to
  full rook placements on a Ferrers board (the permutations correspond
  to square boards), and the construction of the map $\phi^*$ is the
  key step in proving the following result. Let $T$ be a set of
  patterns starting with the prefix $12\cdots k$.  Let $T'$ be the set
  of patterns obtained by replacing this prefix by $k\cdots 21$ in
  every pattern of $T$. Then for all $n$, the number of permutations
  of the symmetric group 
$\Sn_n$ 
  that avoid $T$ equals the number of permutations of $\Sn_n$ that
  avoid $T'$.
  
  Our commutation result, generalized to Ferrers boards, implies that
  the number of {\em involutions\/} of $\Sn_n$ that avoid $T$ is equal
  to the number of involutions of $\Sn_n$ avoiding $T'$, as recently
  conjectured by Jaggard.
\end{abstract}

\maketitle \thispagestyle{empty}

%
%

\section{Introduction}
Let $\pi=\pi_1 \pi_2 \cdots \pi_n$ be a \p\ of length $n$. Let $\tau=
\tau_1 \cdots \tau _k$ be another \p .  An \emm occurrence, of $\tau$
in $\pi$ is a subsequence $\pi_{i_1}\cdots \pi_{i_k}$ of $\pi$ that is
order-isomorphic to $\tau$. For instance, $2 4 6$ is an occurrence of
$\tau= 1 2 3$ in $\pi=2 5 1 43 6$. We say that $\pi$ \emm avoids,
$\tau$ if $\pi$ contains no occurrence of $\tau$. For instance, the
above \p\ $\pi$ avoids $1 2 3 4$. The set of \ps\ of length $n$ is
denoted by $\Sn_n$, and $\Sn_n(\tau)$ denotes the set of
$\tau$-avoiding \ps\ of length $n$.

The idea of systematically studying pattern avoidance in \ps\ appeared
in the mid-eighties~\cite{rodica}. The main problem in this field is
to determine $S_n(\tau)$, the cardinality of $\Sn_n(\tau)$, for any
given pattern $\tau$. This question has subsequently been generalized
and refined in various ways (see for
instance~\cite{babson-einar,italiens,claesson-mansour,krattenthaler},
and~\cite{kitaev-mansour} for a recent survey). However, relatively
little is known about the original question. The case of patterns of
length $4$ is not yet completed, since the pattern $1324$ still
remains unsolved.
See~\cite{bona,gessel-symmetric,stankova,stankova-4,west-these} for
other patterns of length $4$. 

For length $5$ and beyond, all the 
solved cases follow from three important generic results. The first
one, due to Gessel~\cite{gessel-symmetric,gessel-weinstein}, gives the
\gf\ of the numbers $S_n(12\cdots k)$. The second one, due to Stankova
and West~\cite{stankova-west}, states that $S_n(231
\tau)=S_n(312\tau)$ for any pattern $\tau$ on $\{4,5, \ldots , k\}$.
The third one, due to Backelin, West and Xin~\cite{bwx}, shows that
$S_n(12\cdots k \tau)=S_n(k\cdots 21\tau)$ for any pattern $\tau $ on
the set $\{k+1, k+2, \ldots , \ell\}$.  In the present paper an
analogous result is established for pattern-avoiding \emm
involutions,. We denote by $\I_n(\tau)$ the set of involutions
avoiding $\tau$, and by $I_n(\tau)$ its cardinality.

The systematic study of pattern avoiding involutions was also
initiated
in~\cite{rodica}, continued
in~\cite{gessel-symmetric,gouyou-tableaux}
for increasing patterns, and then by Guibert in his
thesis~\cite{guibert-these}. Guibert discovered experimentally that,
for a surprisingly large number of patterns~$\tau$ of length $4$,
$I_n(\tau)$ is the $n$th Motzkin number:
$$
M_n =\sum_{k=0}^{\lfloor n/2\rfloor} \frac{n!}{k!(k+1)!(n-2k)!}.
$$
This was already known for $\tau=1234$ (see \cite{regev}), and
consequently for $\tau=4321$, thanks to the properties of the
Schensted correspondence~\cite{schensted}. Guibert 
explained all the other instances of the Motzkin numbers, except for
two of them: $2143$ and $3214$. However, he was able to describe a
two-label generating tree for the class $\I_n(2143)$. Several years
later, the Motzkin result for the pattern $2143$ was at last derived
from this tree: first in a bijective way~\cite{guibert-vexillaires},
then using \gfs~\cite{mbm-motifs}. No simple generating tree could be
described for involutions avoiding $3214$, and it was only in 2003
that Jaggard~\cite{jaggard} gave a proof of this final conjecture,
inspired by~\cite{babson-west}.  
More generally, he proved that for $k=2$ or 3, $I_n(12\cdots k
\tau)=I_n(k\cdots 21\tau)$ for all 
$\tau$. He conjectured that this holds for all $k$,  which we prove here.

\smallskip We derive this from another result, which may be more
interesting than its implication in terms of forbidden patterns. This
result deals with a transformation $\phi^*$ that was defined
in~\cite{bwx} to prove that $S_n(12\cdots k \tau)=S_n(k\cdots
21\tau)$. This transformation acts not only on \ps , but on more
general objects called \emm full rook placements on a Ferrers shape,
(see Section~\ref{section-wilf} for precise definitions). The map
$\phi^*$ may, at first sight, appear as an \emm ad hoc, construction,
but we prove that it has a remarkable, and far from obvious, property:
it commutes with the inversion of a \p, and more generally with the
corresponding diagonal reflection of a full rook placement.  (By the
inversion of a permutation $\pi$ we mean the map that sends $\pi$,
seen as a bijection, to its inverse.)

The map $\phi^*$ is defined by iterating a transformation $\phi$,
which chooses a certain occurrence of the pattern $k \cdots 21$ and replaces
it by an occurrence of $(k-1) \cdots 21 k$. The map $\phi$ itself does
\emm not, commute with the inversion of \ps , and our proof of the
commutation theorem is actually quite complicated.

This strongly suggests that we need a better description of the map
$\phi^*$, on which the commutation theorem would become obvious. By
analogy, let us recall what happened for the Schensted correspondence:
the fact that the inversion of \ps \ exchanges the two tableaux only
became completely clear with Viennot's description of the
correspondence~\cite{viennot}.

Actually, since the Schensted correspondence has nice properties
regarding the monotone subsequences of \ps , and provides one of the
best proofs of the identity $I_n(12\cdots k)=I_n(k \cdots 21)$, we 
suspect that the map $\phi^*$ might be related to this
correspondence, or to an extension of it to rook placements.

\section{Wilf equivalence for involutions}
\label{section-wilf}
\begin{figure}[tb]
  \centerline{\epsfbox{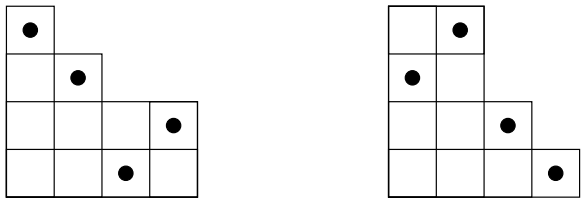}}
\caption{A full rook placement on a Ferrers board, and its inverse.}
\label{fig-placement}
\end{figure}

One of the main implications of this paper is the following.
\begin{thm}\label{thm-involutions}
  Let $k\ge 1$. Let $T$ be a set of patterns, each starting with the
  prefix $12\cdots k$.  Let $T'$ be the set of patterns obtained by
  replacing this prefix by $k\cdots 21$ in every pattern of $T$. Then,
  for all $n \ge 0$, the number of involutions of $\Sn_n$ that avoid
  $T$ equals the number of involutions of $\Sn_n$ that avoid $T'$.
  
  In particular, the involutions avoiding $12\cdots k \tau$ and the
  involutions avoiding $k\cdots 21 \tau$ are equinumerous, for any
  permutation $\tau$ of $\{k+1, k+2, \ldots , \ell\}$.
\end{thm}
This theorem was proved by Jaggard for $k=2$ and $k=3$~\cite{jaggard}.
It is the analogue, for involutions, of a result recently
proved by Backelin, West and Xin for permutations~\cite{bwx}.  Thus it
is not very suprising that we follow their approach. This approach
requires looking at pattern avoidance for slightly more general
objects than permutations, namely, \emm full rook placements on a
Ferrers board., 

Let $\lambda$ be an integer partition, which we represent as a Ferrers
board (Figure~\ref{fig-placement}). A full rook placement on
$\lambda$, or a \emm placement, for short, is a distribution of dots
on this board, such that every row and column contains exactly one
dot. This implies that the board has as many rows as columns.

Each cell of the board will be denoted by its coordinates: in the
first placement of Figure~\ref{fig-placement}, there is a dot in the
cell $({1,4})$. If the placement has $n$ dots, we associate with it a
\p\ $\pi$ of $\Sn_n$, defined by $\pi(i)=j$ if there is a dot in the
cell $(i,j)$. The permutation corresponding to the first placement of
Figure~\ref{fig-placement} is $\pi=4312$.  This induces a bijection
between placements on the $n\times n$ square and permutations of
$\Sn_n$.

The \emm inverse, of a placement $p$ on the board $\lambda$ is the
placement $p'$ obtained by reflecting $p$ and $\lambda$ with respect
to the main diagonal; it is thus a placement on the \emm conjugate, of
$\lambda$, usually denoted by $\lambda'$. This terminology is of course
an extension to placements of the classical terminology for \p s.

\begin{defin}
Let $p$ be a placement on the board $\lambda$, and let $\pi$ be the
corresponding \p. Let $\tau$ be a \p\ of $\Sn_k$. We say that $p$
{\em contains} $\tau$ if 
there exists in $\pi$ an occurrence
$\pi_{i_1}\pi_{i_2}\cdots \pi_{i_k}$ of $\tau$ such that the
corresponding dots are contained in a \emm rectangular, sub-board of
$\lambda$. In other words, the cell with coordinates $(i_k,
\max_j\pi_{i_j})$ must belong to $\lambda$. 
\end{defin}
The placement of
Figure~\ref{fig-placement} contains the pattern $12$, but avoids the
pattern $21$, even though the associated permutation $\pi=4312$
contains several occurrences of $21$. We denote by $\Sn_\lambda(\tau)$
the set of placements on $\lambda$ that avoid $\tau$. If $\lambda$ is
self-conjugate,
we denote by $\I_\lambda(\tau)$ the set of \emm symmetric, (that is,
self-inverse) placements on $\lambda$ that avoid $\tau$. We denote by
$S_\lambda(\tau)$ and $I_\lambda(\tau)$ the cardinalities of these
sets.

In \cite{babson-west,bwx,stankova-west}, it was shown that the notion
of pattern avoidance in 
placements is well suited to deal with prefix exchanges in
 patterns. This was adapted by Jaggard~\cite{jaggard} to
involutions:
\begin{prop}\label{prefix-involutions}
  Let $\alpha$ and $\beta$ be two involutions of $\Sn_k$. Let
  $T_\alpha$ be a set of patterns, each beginning with $\alpha$. Let
  $T_\beta$ be obtained by replacing, in each pattern of $T_\alpha$,
  the prefix $\alpha$ by $\beta$. 
 If, for every self-conjugate shape
  $\lambda$, $I_\lambda(\alpha)= I_\lambda(\beta)$, then
  $I_\lambda(T_\alpha)= I_\lambda(T_\beta)$ for every self-conjugate
  shape.
\end{prop}
Hence Theorem~\ref{thm-involutions} will be proved if we can prove
that $I_\lambda(12\cdots k)= I_\lambda(k\cdots 21)$ for any
self-conjugate shape $\lambda$. A simple induction on $k$, combined
with Proposition~\ref{prefix-involutions}, shows that it is actually
enough to prove the following:
\begin{thm}\label{thm-placements}
  Let $\lambda$ be a self-conjugate shape. Then $I_\lambda(k\cdots
  21)= I_\lambda((k-1) \cdots 21k)$.
\end{thm}
A similar result was proved in~\cite{bwx} for general (asymmetric)
placements: for every shape $\lambda$, one has $S_\lambda(k\cdots 21)=
S_\lambda((k-1) \cdots 21k)$. The proof relies on the description of a
recursive bijection between the sets $\Sn_\lambda(k\cdots 21)$ and
$\Sn_\lambda((k-1) \cdots 21k)$. What we prove here is that this
complicated bijection actually \emm commutes with the inversion
of a placement,, and this implies
Theorem~\ref{thm-placements}.

But let us first describe (and slightly generalize) the transformation
defined by Backelin, West and Xin~\cite{bwx}. This transformation
depends on $k$, which from now on is supposed to be fixed. Since
Theorem~\ref{thm-placements} is trivial for $k=1$, we assume $k \ge2$.
\begin{defin}[\bf{The transformation $\phi$}]\label{dfn-phi}
  Let $p$ be a placement containing $k \cdots 21$, and let $\pi$ be
  the associated \p . To each occurrence of $k \cdots 21$ in $p$,
  there corresponds a decreasing subsequence of length $k$ in $\pi$.
  The \emm \aseq,\ of $p$, denoted by $\A(p)$, is the smallest of these
  subsequences for the lexicographic order.
  
  The corresponding dots in $p$ form
  an occurrence of $k \cdots 21$. Rearrange these dots cyclically
so as to form an
  occurrence of $(k-1) \cdots 21k$. The resulting placement is defined
  to be $\phi(p)$.
  
  If $p$ avoids $k \cdots 21$, we simply define $\phi(p):=p$.  The
  transformation $\phi$ is also called the \emm \A-shift.,
\end{defin}
\begin{figure}[htb]
\begin{center}
  \input{phi-def.pstex_t}
\end{center}
\caption{The \A-shift on the \p\ $7\ 4\ 6\ 3\ 5\ 2\ 1$, when $k=4$.}
\label{fig:phi-def}
\end{figure}

An example is provided by Figure~\ref{fig:phi-def} (the letters of the
\aseq \ are underlined, 
and the corresponding dots are black). It is easy to see that the
\A-shift decreases 
the inversion number of the \p\ associated with the placement 
(details will be given in the proof of Corollary~\ref{confluence}). This
implies that after finitely many iterations of $\bwx$, there will
be no more decreasing subsequences of length $k$ in the placement. We
denote by $\phi ^*$ the iterated transformation, that recursively
transforms \emm every, pattern $k\cdots21$ into $(k-1)\cdots21k$.
For instance, with the \p\ $\pi=7\ 4\ 6\ 3\ 5\ 2\ 1$ of
Figure~\ref{fig:phi-def} and $k=4$, we find
$$
\pi=7\ \un 4\ 6\ \un 3\ 5\ \un 2 \ \un 1 \longrightarrow \un 7\ \un 3\ 
6\ \un 2\ 5\ \un 1\ 4\longrightarrow 3\ 2\ 6\ 1\ 5\ 7\ 4 =\phi^*(\pi).
$$
The main property of $\phi ^*$ that was proved and used in~\cite{bwx} is the
following:
\begin{thm}[\bf{The BWX bijection}]\label{thm-bwx}
For every shape $\lambda$,   the transformation $\phi ^*$
 induces a bijection from
  $\Sn_\lambda((k-1)\cdots 21k )$ to $\Sn_\lambda(k\cdots 21)$.
\end{thm}
\noindent 
The key to our paper is the following rather unexpected theorem.
\begin{thm}[\bf{Global commutation}]\label{thm-global}
  The transformation $\phi ^*$ commutes with the inversion of a
  placement.
\end{thm}
For instance, with $\pi$ as above, we have
$$
\pi^{-1}=7\ \un 6\ \un 4\ \un 2\ 5\ 3\ \un 1 \longrightarrow \un 7\ \un
4\ \un 2\ \un 1\ 5\ 3\ 6 \longrightarrow 4 \ 2\ 1\ 7\ 5\ 3\ 
6=\phi^*(\pi^{-1})
$$
and we observe that
$$
\phi^*(\pi^{-1})=\left(\phi^*(\pi)\right)^{-1}.
$$
Note, however, that $\phi(\pi^{-1})\not
=\left(\phi(\pi)\right)^{-1}$. Indeed, $\phi(\pi^{-1})= 7\ 4\ 2\ 1\ 5\ 
3\ 6 $ while $\left(\phi(\pi)\right)^{-1} =6\ 4\ 2\ 7\ 5\ 3\ 1$, so
that the elementary transformation $\phi$, that is, the \A-shift, does
not commute with the inversion.

Theorems~\ref{thm-bwx} and~\ref{thm-global} together imply that
$\phi^*$ induces a bijection from\lb $\I_\lambda((k-1)\cdots 21k )$ to
$\I_\lambda(k\cdots 21)$, 
for every self-conjugate shape $\lambda$. This proves Theorem~\ref{thm-placements},
and hence Theorem~\ref{thm-involutions}. The rest of the paper is
devoted to proving Theorem~\ref{thm-global}, which we call the theorem
of \emm global commutation., By this, we mean that the inversion
commutes with the global tranformation $\phi ^*$ (but not with the
elementary transformation $\phi $).

\medskip
\noindent{\bf Remarks}\\
1. At first sight, our definition of the \aseq\ 
(Definition~\ref{dfn-phi}), does not seem to coincide with the
definition given in~\cite{bwx}. Let $a_k \cdots a_2 a_1$ denote the
\aseq\ of the placement $p$, with $a_k > \cdots > a_1$. We identify
this sequence with the corresponding set 
of dots in $p$.  The dot $a_k$ is the lowest dot that is the leftmost
point in an occurrence of $k\cdots 21$ in $p$.  Then $a_{k-1}$ is the
 lowest dot such that $a_ka_{k-1}$ is the beginning of an
occurrence of $k\cdots 21$ in $p$, and so on.

However, in \cite{bwx}, the dot $a_k$ is chosen as above, but then
each of the next dots $a'_{k-1}, \ldots , a'_1$ is chosen to be as far
\emm left, as possible, and not as \emm low, as possible. Let us prove
that the two procedures give the same sequence of dots. Assume not,
and let $a_j\not = a'_j$ be the first (leftmost) point where the two
sequences differ. By definition, $a_j$ is lower than $a'_j$, and to
the right of it. But then the sequence $a_{k-1} \cdots a_{j+1}a'_j a_j
\cdots a_2a_1$ is an occurrence of the pattern $k \cdots 21$ in $p$,
which is smaller than $a_k \cdots a_2 a_1$ for the lexicographic
order, a contradiction.

The fact that the \aseq\ can be defined in two different ways will be
used very often in the paper.

\noindent 2. At this stage, we have reduced the proof of
Theorem~\ref{thm-involutions} to the proof of the global commutation
theorem, Theorem~\ref{thm-global}.

\section{From local commutation to global commutation}
In order to prove that $\phi^*$ commutes with the inversion of
placements, it would naturally be tempting to prove that $\phi$ itself
commutes with the inversion. However, this is not the case, as shown
above. Given a placement $p$ and its inverse $p'$, we thus want to
know how the placements $\phi(p)$ and $\phi(p')'$ differ.

\begin{defin}
  For any shape $\lambda$ and any placement $p$ on $\lambda$, we
  define $\BWX(p)$ by
  $$
  \BWX(p):= \phi(p')'.
  $$
  Thus $\BWX(p)$ is also a placement on $\lambda$.
\end{defin}
Note that $\BWX^m(p)=(\bwx^m(p'))'$, so that the theorem of global
commutation, Theorem~\ref{thm-global}, can be restated as
$\BWX^*=\bwx^*$.

Combining the above definition of $\BWX$ with Definition~\ref{dfn-phi}
gives an alternative description of $\BWX$.
\begin{lemma}[\bf{The transformation $\psi$}]\label{dfn-psi}
  Let $p$ be a placement containing $k \cdots 21$.
Let $b_1,b_2, \ldots ,  b_k$ be defined recursively as follows:  
For all $j$, $b_j$ is the leftmost dot
  such that $b_j \cdots b_2 b_1$ ends an occurrence of $k \cdots 21$
  in $p$. We call $b_k \cdots b_2 b_1$ the \emm\bseq, of $p$, and denote
  it by $\B(p)$.
  
  Rearrange the $k$ dots of the \bseq\ cyclically so as to form an
  occurrence of\lb $(k-1) \cdots 21k$: the resulting placement is
  $\psi(p)$.
  
  If $p$ avoids $k \cdots 21$, then $\psi(p)=p$.  The transformation
  $\psi$ is also called the \emm\B-shift,.
\end{lemma}
According to the first remark that concludes Section~\ref{section-wilf}, we
can alternatively define $b_j$, for $j\ge 2$, as the \emm lowest, dot
such that $b_j \cdots b_2 b_1$ ends an occurrence of $k \cdots 21$ in
$p$.

We have seen that, in general, $\phi$ does not commute with the
inversion. That is, $\phi(p)\not = \BWX(p)$ in general. The above
lemma tells us that $\phi(p) = \BWX(p)$ if and only if the \aseq\ and
the \bseq\ of $p$ coincide. If they do \emm not, coincide, then we
still have the following remarkable property, whose proof is deferred
to the very end of the paper.
\begin{thm}[\bf{Local commutation}]\label{thm-local}
  Let $p$ be a placement for which the \A- and \bseq s do not
  coincide. 
Then $\phi(p)$ and $\psi(p)$ still contain the pattern $k\cdots 21$, and 
  $$
  \bwx(\BWX(p))=\BWX(\bwx(p)).
  $$
\end{thm}
For instance, for the \p\ of Figure~\ref{fig:phi-def} and $k=4$, we have the
following commutative diagram, in which the underlined
(resp.~overlined) letters correspond to the \aseq\ (resp.~\bseq):\\
\begin{center}
  \input{local.pstex_t}
\end{center}
\smallskip

A classical argument, which is sometimes stated in terms of \emm
locally confluent, and \emm globally confluent rewriting systems,
(see~\cite{huet} and references therein), will show that
Theorem~\ref{thm-local} implies $\BWX^*=\bwx^*$, and actually the more
general following corollary.
\begin{cor}\label{confluence}
  Let $p$ be a placement. Any iterated application of the
  transformations $\bwx$ and $\BWX$ yields ultimately the same
  placement, namely $\phi^*(p)$. Moreover, all the minimal sequences
  of transformations that yield $\phi^*(p)$ have the same length.
\end{cor}
%
\begin{figure}[b]
\hspace{-8mm}
  \input{component-final.pstex_t}
\caption{The action of  $\bwx$ and $\BWX$ on a part of $\Sn_9$, for $k=4$.}
\label{fig-connected}
\end{figure}

Before we prove this corollary, let us illustrate it. We think of the
set of \p s of length $n$ as the set of vertices of an oriented graph,
the edges of which are given by the maps $\bwx$ and $\BWX$.
Figure~\ref{fig-connected} shows a connected component of this graph.
The dotted edges represent $\bwx$ while the plain edges represent
$\BWX$. The dashed edges correspond to the cases where $\bwx$ and
$\BWX$ coincide.  We see that all the paths that start at a given
point converge to the same point.

\begin{proof}
  For any placement $p$, define the inversion number of $p$ as the
  inversion number of the associated \p \ $\pi$ (that is, the number
  of pairs $(i,j)$ such that $i<j$ and $\pi_i>\pi_j$). Assume $p$
  contains at least one occurrence of $k\cdots 21$, and let $i_1 <
  \cdots <i_k$ be the positions (abscissae) of the elements of the
  \aseq\ of $p$. A careful examination of the inversions of $p$ and
  $\bwx(p)$ shows that
  $$
  \inv(p)-\inv(\bwx(p)) = k-1 + 2 \sum_{m=1}^{k-1} \hbox{Card }\{ i:
  i_m<i<i_{m+1} \mbox{~and~} \ \pi_{i_1} >\pi_i > \pi_{i_{m+1}} \}.
  $$
  In particular, $\inv(\bwx(p))<\inv(p)$. By symmetry, together with the
  fact that $\inv(\pi^{-1})= \inv(\pi)$, it follows that
  $\inv(\BWX(p))<\inv(p)$ too.
  
  We encode the compositions of the maps $\bwx$ and $\BWX$ by words on
  the alphabet $\{\phi, \psi\}$. For instance, if $u$ is the word
  $\bwx\BWX^2$, then $u(p)=\bwx\BWX^2(p)$.
Let us prove, by induction on $\inv(p)$, the following two
  statements:\\
  1. If $u$ and $v$ are two words such that $u(p)$ and $v(p)$ avoid
  $k\cdots 21$,
  then  $u(p)=v(p)$. \\
  2. Moreover, if $u$ and $v$ are minimal for this property (that is,
   for any non-trivial factorization $u=u_0u_1$, the placement
  $u_1(p)$ still contains an occurrence of $k\cdots 21$ --- and
  similarly for $v$), then $u$ and $v$ have the same length.

If the first property holds for $p$, then $u(p)=v(p)= \bwx^*(p)$. If
  the second property holds, we denote  by $L(p)$ the length of any
  minimal word $u$ such that $u(p)$ avoids $k\cdots21$.

  If $\pi$ is the identity, then the two results are obvious. They
  remain obvious, with $L(p)=0$, if $p$ does not contain any
  occurrence of $k\cdots 21$.
  
  Now assume $p$ contains such an occurrence, and $u(p)$ and $v(p)$
  avoid $k\cdots 21$. By assumption, neither $u$ nor $v$ is the empty
  word. Let $f$ (resp.~$g$) be the rightmost letter of $u$
  (resp.~$v$), that is, the \emm first, transformation that is applied
  to $p$ in the evaluation of $u(p)$ (resp. $v(p)$). Write $u=u'f$ and
  $v=v'g$.
  
  If $f(p)=g(p)$, let $q$ be the placement $f(p)$. Given that
  $\inv(q)<\inv(p)$, and that the placements $u(p)=u'(q)$ and $v(p)=v'(q)$
  avoid $k\cdots 21$, both statements follow by induction. 
  
  If $f(p)\not =g(p)$, we may assume, without loss of generality, that
  $f=\bwx$ and $g=\BWX$. Let $q_1=\bwx(p)$, $ q_2=\BWX(p)$ and
  $q=\bwx(\BWX(p)=\BWX(\bwx(p))$ (Theorem~\ref{thm-local}).  The
  induction hypothesis, applied to $q_1$, gives
  $u'(q_1)=\phi^*(\psi(q_1))=\phi^*(q)$, that is, $u(p)=\phi^*(q)$
  (see the figure below).
  Similarly, $v'(q_2)=\phi^*(q_2)=\phi^*(q)$, that is,
  $v(p)=\phi^*(q)$. This proves the first statement. If $u$ and $v$
  are minimal for $p$, then so are $u'$ and $v'$ for $q_1$ and $q_2$
  respectively. By the first statement of Theorem~\ref{thm-local}, $q_1$ and
  $q_2$ still contain the pattern $k\cdots 21$, so $L(q)=L(q_1)-1=
  L(q_2)-1$, and  the words $u'$ 
  and $v'$ have the same length. Consequently, $u$ and $v$ have the
  same length too.
\end{proof}

\begin{center}
  \input{local-global.pstex_t}
\end{center}

\noindent
{\bf Note.} We have reduced the proof of Theorem~\ref{thm-involutions}
to the proof of the local commutation theorem,
Theorem~\ref{thm-local}.  The last two sections of the paper are
devoted to this proof, which turns out to be unexpectedly complicated.
There is no question that one needs to find a more illuminating
description of $\phi^*$, or of $\phi \circ \psi$, which makes
Theorems~\ref{thm-global} and~\ref{thm-local} clear.

\section{The local commutation for permutations}
\begin{figure}[htb]
\begin{center}
  \input{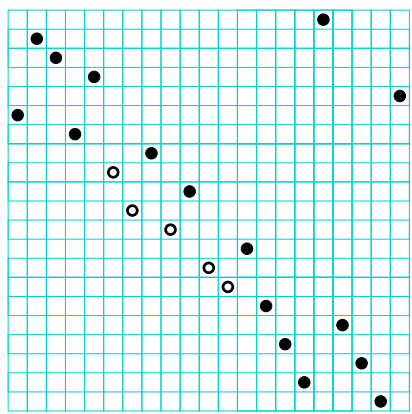}
\end{center}
\caption{Top: A permutation $\pi$, with its \A- and \B-sequences
  shown. Left: After the \B-shift. Right: After the \A-shift. Bottom:
  After the composition of $\phi$ and $\psi$.}
\label{fig-example}
\end{figure}

In this section, we prove that the local commutation theorem holds for
permutations. It will be extended to placements in the next section.
To begin with, let us study a big example, and use it to describe the
contents and the structure of this section. This example is
illustrated in Figure~\ref{fig-example}.

\medskip
\noindent{\bf Example.}
Let $\pi$ be the following permutation of length $21$:
$$
\pi= 17\ 21\ 20\ 16\ 19\ 18\ 13\ 15\ 11\ 14\ 12\ 8\ 10\ 9\ 7\ 4\ 2\ 
6\ 5\ 3\ 1.
$$
1. Let $k=12$. The \aseq \ of $\pi$ is
$$
\A(\pi)= 17\ 16/15\ 14\ 12\ 10\ 9\ 7/6\ 5\ 3\ 1 ,
$$
while its \bseq\ is
$$
\B(\pi)= 21\ 20\ 19\ 18/15\ 14\ 12\ 10\ 9\ 7/4\ 2.
$$
Observe that the intersection of $\A(\pi)$ and $\B(\pi)$ (delimited by
'/') consists of the 
letters $15\ 14\ 12\ 10\ 9\ 7$, and that they are consecutive both in
$\A(\pi)$ and $\B(\pi)$.  Also, \B\ contains more letters than \A\ 
before this intersection, while \A\ contains more letters than \B\ 
after the intersection. We prove that this is always true in
Section~\ref{section-AB} below (Propositions~\ref{contiguous}
and~\ref{a-more-after}).  

\smallskip \noindent 2. Let us now apply the \B-shift to $\pi$. One
finds:
$$\psi(\pi)= 17\ 20\ 19\ 16\ 18\ 15\ 13\ 14\ 11\ 12\ 10\ 8\ 9\ 7\ 4\ 
2\ 21\ 6\ 5\ 3\ 1.
$$
The new \aseq\ is now $\A(\psi(\pi))= 17\ 16 /15\ 13\ 11\ 10\ 8\ 
7/6\ 5\ 3\ 1$. Observe that all the letters of $\A(\pi)$ that were
before or after the intersection with $\B(\pi)$ are still in the new
\aseq , as well as the first letter of the intersection.  We prove
that this is always true in Section~\ref{section-a-after-B}
(Propositions~\ref{a-before-int} and~\ref{a-after-int}).  In this
example, the last letter of the intersection is still in the new \aseq
, but this is not true in general.

By symmetry with respect to the main diagonal, after the \A-shift, the
letters of \B \ that were before or after the intersection are in the
new \bseq , as well as the first letter of \A \ following the
intersection (Corollary~\ref{b-begin-end}). This can be checked on our example:
$$
\phi(\pi)= 
16\ 21\ 20\  15\  19\  18\  13\  14\  11\  12\  10\  8\  9\  7\  6\
4\  2\  5\  3\  1\  17,
%
$$
and the new \bseq\ is $\B(\phi(\pi))= 21\ 20\ 19\ 18/ 13\ 11\ 10\ 
8\ 7\ 6\ / 4\ 2.$

\smallskip \noindent 3. Let $a_i=b_j$ denote the first (leftmost)
point in $\A(\pi)\cap \B(\pi)$, and let $a_d=b_e$ be the last point in this
intersection. We have seen that after the \B-shift, the new \aseq\ 
begins with $a_k\cdots a_i= 17\ 16 \ 15$, and ends with $a_{d-1}\cdots
a_1=6\ 5\ 3\ 1$. The letters in the center of the new \aseq, that is,
the letters replacing $a_{i-1}\cdots a_d$,  are
$x_{i-1}\cdots x_d= 13\ 11\ 10\ 8\ 7$.  Similarly, after the \A-shift,
the new \bseq\ begins with $b_k\cdots b_{j+1}= 21\ 20\ 19\ 18$, and
ends with $a_{d-1}b_{e-1}\cdots b_1= 6\ 4\ 2$. The central letters are
again $x_{i-1}\cdots x_d= 13\ 11\ 10\ 8\ 7$! (See
Figure~\ref{fig-example}). This is not a coincidence; we prove in
Section~\ref{section-composition} that this always holds
(Proposition~\ref{same-int-lemma}). 

\smallskip \noindent 4. We finally combine all these properties to
describe explicitly how the maps $\phi\circ\psi$ and $\psi\circ\phi$
act on a permutation $\pi$, and conclude that they yield the same
permutation if the \A- and \B-sequences of $\pi$ do not coincide
(Theorem~\ref{thm-permutations}). 

\newpage

\subsection{The \aseq \ and the \bseq}\label{section-AB}
\begin{defin}[\bf{Labels}]\label{Labels}
  Let $\pi=\pi_1\pi_2\cdots \pi_n$ be a permutation. For $1\le i \le
  n$, let $\ell_i$ be the maximal length of a decreasing subsequence
  in $\pi$ that starts at $\pi_i$. The \emm length sequence,, or \emm
  $\ell$-sequence,, of $\pi$ is $\ell(\pi)=\ell_1\ell_2\cdots \ell_n$.
  Alternatively, it can be defined recursively as follows: $\ell_n=1$
  and, for $i<n$,
  \begin{equation}\label{label-def}
  \ell_i= \max\{\ell_m\}+1,
  \end{equation}
  where the maximum is taken over all $m>i$ such that $\pi_m <\pi_i$.
  
  We refer to the entries of the \lseq\ as \emm labels, and say that
  the label $\ell_i$ is \emm associated to, the letter $\pi_i$ in
  $\pi$.  Also, if $x=\pi_i$ then, abusing notation, we let
  $\ell(x)=\ell_i$.
  
  Given a subsequence $s = \pi_{i_1}\pi_{i_2}\cdots\pi_{i_k}$ of
  $\pi$,  we say that
  $\ell_{i_1},\ell_{i_2},\ldots,\ell_{i_k}$ is the subsequence of
  $\ell(\pi)$ \emm associated to $s$,.
\end{defin}

Here is an example, where we have written the label of $\pi_i$ below
$\pi_i$ for each $i$:
$$
\begin{array}[c]{cccccccccc}
\pi= & 3 & 7 & 4 & 9 & 1 & 8 & 5 & 6 & 2\\
& 2 & 3 & 2 & 4 & 1 & 3 & 2 & 2 & 1
\end{array}
$$
The subsequence of $\ell(\pi)$ associated to $741$ is $3,2,1$.

\begin{lemma}
  The subsequence of $\ell(\pi)$ associated to a decreasing
  subsequence $x_m \cdots x_2 x_1$ in $\pi$ is strictly decreasing.  In
  particular, $\ell (x_i) \ge i$ for all $i$.
\end{lemma}
\begin{proof}
  Obvious, by definition of the labels.
\end{proof}

\begin{lemma}\label{same-label-increasing}
  Let $x_1,\ldots,x_i$ be, from left to right, the list of letters in
  $\pi$ that have label $m$.  Then $x_1<x_2<\cdots<x_i$.
\end{lemma}
\begin{proof}
  If $x_j>x_{j+1}$ then, since $x_j$ precedes $x_{j+1}$, we would have
  $\ell(x_j)>\ell(x_{j+1})$, contrary to assumption.
\end{proof}

\begin{defin}[\bf{Successor sequence}]
  Let $x$ be a letter in $\pi$, with $\ell(x)=m$.  The {\em successor
    sequence $s_m s_{m-1}\cdots s_1$ of $x$} is the sequence of
  letters of $\pi$ such that $s_m=x$ and, for $i\le m$,
 $s_{i-1}$ is the
  first (leftmost) letter after $s_i$ such that
  $\ell(s_{i-1})=\ell(s_i)-1$.  In this case, we say that $s_{m-1}$ is
  the \emm label successor of $x$.,
\end{defin}

\begin{lemma} The successor sequence of $x$ is a decreasing sequence.
\end{lemma}
\begin{proof}
  By definition of the labels, one of the letters labelled $\ell(x)-1$
  that are to the right of $x$ is smaller than $x$. By
  Lemma~\ref{same-label-increasing}, the leftmost of them, that is,
  the label successor of $x$, is smaller than $x$.
\end{proof}

Let us now rephrase, in terms of permutations, the definitions of the
\aseq \ and \bseq \ (Definition~\ref{dfn-phi} and
Lemma~\ref{dfn-psi}).

Given a permutation $\pi$ that contains a decreasing subsequence of
length $k$, the \aseq\ of $\pi$ is the sequence $\A(\pi)=a_k
a_{k-1} \cdots a_1$, where 
for all $i$, $a_i$ is the \emm smallest, letter in $\pi$ such that $a_k \cdots
a_{i+1}a_i$ is the \emm prefix, of a decreasing sequence of $\pi$ of
length $k$. 
The  \bseq\  of $\pi$ is $\B(\pi)=b_k b_{k-1} \cdots b_1$, where
for all $i$, $b_i$ is the
\emm leftmost, letter such that  $b_i b_{i-1}\cdots b_1$ is the \emm
suffix, of a decreasing sequence of length $k$.  
According to the remark at the end of Section~\ref{section-wilf}, the
letter $a_i$ can alternatively be chosen as \emm left, as possible
(for $i<k$), and the letter $b_i$ as \emm small, as possible (for
$i>1$).

The three simple lemmas above, as well as Lemma~\ref{a-succ}
  below, will be used frequently, but without specific mention, in the
  remainder of this section.  From now on we denote by $\A = a_k
  \cdots a_1$ and $\B = b_k \cdots b_1$ the \A- and \B-sequences of
  $\pi$.  The next lemma characterizes the \aseq \ in terms of labels.
\begin{lemma}
\label{a-succ}
  The letter $a_k$ is the leftmost letter in~$\pi$ with label $k$ and,
  for $i<k$, $a_i$ is the first letter after $a_{i+1}$ that has label
  $i$.  In particular, the \aseq\ of $\pi$ is the successor sequence of $a_k$,
  and the subsequence of $\ell(\pi)$ associated to $\A(\pi)$ is
  $k,k-1,\ldots,1$.
\end{lemma}
\begin{proof}
  Clearly, the label of $a_k$ must be at least $k$. If it is larger than
  $k$, then the label successor of $a_k$ is smaller than $a_k$ and is
  the first letter of a decreasing sequence of length $k$, a
  contradiction. Hence the label of $a_k$ must be exactly $k$. Now
  given that $a_k$ has to be as small as possible,
  Lemma~\ref{same-label-increasing} implies that $a_k$ is the leftmost
  letter having label $k$.
  
  We then proceed by decreasing induction on $i$.  Since $a_i$ is
  smaller than, and to the right of, $a_{i+1}$, its label must be at
  most $i$. Since $a_k\cdots a_i$ is the prefix of a decreasing
  sequence of length $k$, the label of $a_i$ must be at least $i$, and
  hence, exactly $i$. Since we want $a_i$ to be as small as possible,
  it has to be the first letter after $a_{i+1}$ with label $i$
  (Lemma~\ref{same-label-increasing}).
\end{proof}

\begin{lemma}[\bf{The key lemma}]\label{key-lemma}
  Let $i \le j$. Suppose $\pi$ contains a decreasing sequence of the
  form $b_{j+1} x_j \cdots x_i$ such that $x_i$ precedes $b_i$.
  Then $\ell(x_i)<\ell(b_i)$.
\end{lemma}
\begin{proof}
  First, observe that, by definition of \B , one actually has
  $x_i<b_i$ (otherwise, the \bseq\ could be extended).  Suppose that
  $\ell(x_i)\ge\ell(b_i)$.  In particular, then, $\ell(x_i)\ge i$. Let us
  write $\ell(x_i)=i+r$, with $r\ge 0$. Let $x_i x_{i-1} \cdots
  x_{-r+1}$ be the successor sequence of $x_i$, so that
  $\ell(x_p)=p+r$ for all $p\le i$.  Now,
  $\ell(x_{i-m})\ge\ell(b_{i-m})$ for all $m\in[0,i-1]$, because the labels of
  the \bseq\ are strictly decreasing.  Thus, if $x_{i-m}<b_{i-m}$ then
  $x_{i-m}$ must precede $b_{i-m}$, for otherwise
  $\ell(b_{i-m})>\ell(x_{i-m})$.
  
  Recall that $x_i<b_i$. Let $m$ be the largest integer with $m<i$
  such that $x_{i-m}<b_{i-m}$, which implies that $x_{i-m}$ precedes
  $b_{i-m}$ (clearly, $m\ge0$).  If $m=i-1$ then $x_{1}<b_1$, so
  $x_{1}$ precedes $b_1$, and thus the sequence
  $$
  b_k \cdots b_{j+1} x_j x_{j-1} \cdots x_{1}
  $$
  is decreasing, has length $k$ and ends to the left of $b_1$,
  which contradicts the definition of the \bseq .
  Thus $m<i-1$.  
Now,  $x_{i-m}<b_{i-m}$ and $x_{i-m}$ precedes $b_{i-m}$,
  but $x_{i-m-1}>b_{i-m-1}$.
Note that $x_{i-m}$ precedes $b_{i-m-1}$ since it precedes  $b_{i-m}$.
  Thus, the sequence
  $$
  b_k \cdots b_{j+1} x_j x_{j-1} \cdots x_{i-m} b_{i-m-1}
  \cdots b_1
  $$
  is decreasing and has length $k$.  Since $x_{i-m}$ precedes
  $b_{i-m}$, the definition of $\B$ implies that $x_{i-m}$ would have
  been chosen instead of $b_{i-m}$ in $\B$.  This is a contradiction,
  so $\ell(x_i)<\ell(b_i)$.
\end{proof}

\begin{lemma}
\label{new-lemma}
Assume the label successor of $b_m$ does not belong to \B. Then no
letter of the successor sequence of $b_m$ belongs to \B, apart from
$b_m$ itself.
\end{lemma}
\begin{proof}
 Let $x$ be the label successor of $b_m$, and let $x_r \cdots x_1$ be the
  successor sequence of $x$, with $x_r=x$.  Assume one of the $x_i$ belongs to \B,
  and let $x_{s-1}=b_j$ be the leftmost of these. The
  successor sequence of $b_m$ thus reads $b_m x_r \cdots x_s
  b_j x_{s-2} \cdots x_1$. By assumption, $s\le r$.
  
  We want to prove that the sequence $x_r \cdots x_s$ is longer than
  $b_{m-1} \cdots b_{j+1}$, which will contradict the definition of
  the \bseq\ of $\pi$. We have $\ell(x_r)+1=\ell(b_m)>\ell(b_{m-1})$,
  so that $ \ell(x_r) \ge \ell(b_{m-1})$. Hence by Lemma~\ref{key-lemma},
  $b_{m-1}$ precedes $x_r$. This implies that
  $$
  \ell(x_r)> \ell(b_{m-1}),
  $$
  for otherwise the label successor of $b_m$ would be $b_{m-1}$ instead
  of $x_r$. At the other end of the sequence $x_r \cdots x_s$, we
  naturally have
  $$
  \ell(x_s)=1+\ell(b_j) \le \ell(b_{j+1}).
  $$
  Given that the $x_i$ form a successor sequence, while the labels
  of \B\ are strictly decreasing, the above two inequalities imply
  that $x_r \cdots x_s$ is longer than $b_{m-1} \cdots b_{j+1}$, as
  desired.
\end{proof}

\begin{lemma}\label{b-after-int}
  Assume that $a_d=b_e$ with $e>1$, and that $b_{e-1}$ does not belong to $\A$.
  Then $d>1$ and $b_{e-1}$ precedes $a_{d-1}$.  Moreover,
  $b_{e-1}<a_{d-1}$ and $\ell(b_{e-1})<\ell(a_{d-1})=d-1$.
  
  By symmetry, if $a_i =b_j$ with $i<k$ and $a_{i+1}$ does not belong
  to \B, then 
  $j<k$ and $a_{i+1}$ precedes $b_{j+1}$. Moreover, $a_{i+1}<b_{j+1}$.
\end{lemma}
\begin{proof}
  If $d=1$, then $a_k \cdots a_1$ is a decreasing sequence of length
  $k$ that ends to the left of $b_1$ and this contradicts the definition
  of \B . Hence $d>1$.
  
  We have $a_{d-1}<b_{e}$ and $\ell(a_{d-1})\ge\ell(b_{e-1})$.  By
  Lemma~\ref{key-lemma}, this implies that $b_{e-1}$ precedes
  $a_{d-1}$.
  
  By the definition of \A, we have $b_{e-1}<a_{d-1}$, for otherwise
  $b_{e-1}$ could be inserted in $\A$.
  
  Finally, if $\ell(b_{e-1})=\ell(a_{d-1})$ then $b_{e-1}$ would be
  the next letter after $a_d$ in the \aseq, since $b_{e-1}$ precedes
  $a_{d-1}$.
\end{proof}

\begin{prop}
\label{contiguous}
  If $b_e\in\A$ and $b_{e-1}\not\in\A$ then $b_m\not\in\A$ for all
  $m<e$.  Consequently, the intersection of \A\ and \B\ is a
  (contiguous) segment of each sequence.
\end{prop}
\begin{proof}
  Suppose not, so there is an $m<e-1$ with $b_m\in\A$.  Let $d$ and
  $p$ be such that $a_d=b_e$ and $a_p=b_m$.  By
  Lemma~\ref{b-after-int}, $\ell(b_{e-1})<\ell(a_{d-1})$, so there are
  more letters in the \aseq\ than in the \bseq\ between $b_e$ and
  $b_m$.  But then the sequence
  $$
  b_k \cdots b_e a_{d-1} \cdots a_p b_{m-1} \cdots b_2
  $$
  has length at least $k$, which contradicts the definition of the
  \bseq.
  Hence the intersection of \A\ and \B\ is formed of consecutive
  letters of \B. By symmetry, it also consists of consecutive letters
  of \A.
\end{proof}

The preceding proposition will be used implicitly in the remainder of this
section.

\begin{prop}\label{a-more-after}
  If \A\ and \B\ intersect but do not coincide, then the \aseq\ contains
  more letters than the \bseq\ after the intersection and the \bseq\ 
  contains more letters than the \aseq\ before the intersection.  In
  particular, if $a_1$ belongs to \B\ or $b_k$ belongs to \A, then the
  \aseq\ and the \bseq\ coincide.
\end{prop}
\begin{proof}
  Let $b_e=a_d$ be the last letter of the intersection.  The \aseq\ 
  has exactly $d-1$ letters after the intersection.  Let us first
  prove that $d\ge2$. Assume $d=1$. Then by Lemma~\ref{b-after-int},
  $e=1$. Let $a_i=b_i$ be the largest element of $\A\cap \B$. By
  assumption, $i<k$. By Lemma~\ref{b-after-int}, $a_{i+1}$ precedes
  $b_{i+1}$ and is smaller. This contradicts the definition of \B .
  Hence $d >1$.
  
  If the \bseq\ contains any letters after the intersection, then
  $\ell(b_{e-1})<d-1$, according to Lemma~\ref{b-after-int}, so the
  \bseq\ can contain at most $d-2$ letters after the intersection.
  This proves the first statement. The second one follows by symmetry
  (or by subtraction).
\end{proof}
 
\begin{lemma}\label{bk-label}
  Assume $\ell(b_k)=k$. Then the \aseq\ and \bseq\ coincide.
\end{lemma}
\begin{proof}
  Assume the two sequences do not coincide.  If they intersect, their
  last common point being $b_e=a_d$, then
  Proposition~\ref{a-more-after} shows that the sequence 
$$
b_k\cdots  b_e a_{d-1}\cdots a_1
$$ 
is decreasing and has length $>k$. This
  implies that $\ell (b_k) > k$, a contradiction.
  
  Let us now assume that the two sequences do not intersect.  By
  definition of the \aseq , $a_k<b_k$. If $b_k$ precedes $a_k$, then
  $\ell(b_k)>\ell(a_k)=k$, another contradiction.
  Thus $a_k$ precedes $b_k$. Let us prove by decreasing
  induction on $h$ that $a_h$ precedes $b_h$ for all $h$. If this is
  true for some $h \in [2, k]$, then $a_h$ is in the \aseq , $a_{h-1}$
  and $b_{h-1}$ lie to its right and have the same label. Since
  $a_{h-1}$ is chosen in the \aseq , it must be left of $b_{h-1}$.  By
  induction, we conclude that $a_1$ precedes $b_1$, which contradicts
  the definition of the \bseq.
\end{proof}

\smallskip
\subsection{The \aseq\ after the $\B$-shift}
\label{section-a-after-B}
 We still denote by
$\A = a_k \cdots a_1$ and $\B = b_k \cdots b_1$ the \A- and
\B-sequences of a permutation $\pi$. Recall that the \B-shift performs a cyclic
shift of the elements of the \bseq , and is denoted $\psi$.
We begin with a sequence of lemmas that tell us how the labels evolve
 during the \B-shift.
\begin{lemma}[\bf{The order of \A}]\label{a-in-order}
  Assume \A \ and \B\ do not coincide. In $\BWX(\pi)$, the letters
  $a_k, \ldots , a_2, a_1$ appear in this order.
In particular, $\ell(a_i)\ge i$ in $\psi(\pi)$, and $\psi(\pi)$
  contains the pattern $k\cdots 21$.
\end{lemma}
\begin{proof}
  The statement is obvious if \A\ and \B\ do not intersect.
  Otherwise, let $a_i=b_j$ be the first (leftmost) letter of $\A\cap \B$ and let
  $a_d=b_e$ be the last letter of $\A\cap \B$.  By
  Proposition~\ref{a-more-after}, $j<k$.  Hence when we do the
  \B-shift, the letters $a_i, \ldots , a_d$ move to the left, while
  the other letters of \A\ do not move.  Moreover, the letter $a_{i}$
  will replace $b_{j+1}$, which, by Lemma~\ref{b-after-int}, is to the
  right of $a_{i+1}$.
  Hence the letters $a_k, \ldots , a_2, a_1$ appear in this order
  after the \B-shift.
\end{proof}

\begin{lemma}\label{incr-label}
  Let $x\le b_k$. Then the label of $x$ cannot be larger in
  $\BWX(\pi)$ than in $\pi$.
\end{lemma}
\begin{proof}
  We proceed by induction on $x \in \{1, \ldots , b_k\}$, and use the
  definition~\Ref{label-def} (in Definition \ref{Labels}) of the labels. The result is obvious for
  $x=1$. Take now $x\ge 2$, and assume the labels of $1,\ldots , x-1$
  have not increased. If $x\not\in\B$, all the letters that are
  smaller than $x$ and to the right of $x$ in $\BWX(\pi)$ were already
  to the right of $x$ in $\pi$, and have not had a label increase by
  the induction hypothesis. Thus the label of $x$ cannot have
  increased. The same argument applies if $x=b_k$.
 
  Assume now that $x=b_m$, with $m<k$. Then $b_m$ has moved to the
  place of $b_{m+1}$ during the \B-shift. The letters that are smaller
  than $b_m$ and were already to the right of $x$ in $\pi$ have not
  had a label increase. Thus they cannot entail a label increase for
  $b_m$.  The letters that are smaller than $b_m$ and lie between
  $b_{m+1}$ and $b_m$ in $\pi$ have label at most $\ell(b_m)-1$ in
  $\pi$ (Lemma~\ref{key-lemma}), and hence in $\psi(\pi)$, by the
  induction hypothesis. Thus they cannot entail a label increase for
  $b_m$ either. Consequently, the label of $b_m$ cannot change.
\end{proof}

Note that the label of letters larger than $b_k$ may increase, as
shown by the following example, where $k=3$:
$$
\begin{array}[c]{cccccccccccccccccccccccccccccccccccccccc}
\pi &=& 3 & \ov 7& \ov 4&  {\bf 8}&  \ov 1 & 5& 6& 2 &
\rightarrow &\psi(\pi)&=& 
3 & 4 & 1 &  {\bf 8} & 7 & 5 & 6 & 2 \\
    & & 2 & 3 & 2 & {\bf 3} & 1 & 2 & 2 & 1 
&&&& 2 &2 &1 &{\bf 4}& 3& 2& 2& 1.
\end{array}
$$

\begin{lemma}[\bf{The labels of \A}]\label{a-labels-no-change}
  Assume \A\ and \B\ do not coincide.  The labels associated to the
  letters $a_k, \ldots, a_1$ do not change during
the $\B$-shift.
\end{lemma}
\begin{proof}
  By Lemma~\ref{a-in-order}, the label of $a_i$ cannot decrease and by
  Lemma~\ref{incr-label}, it cannot increase either.
\end{proof}

\begin{lemma}[\bf{The labels of \B}]\label{no-b-change}
  Let $m<k$. The label associated to $b_m$ does not change during the
  \B-shift (although $b_m$ moves left).
\end{lemma}
\begin{proof}
  By Lemma~\ref{incr-label}, the label of $b_m$ cannot increase.  Assume
  that it decreases, and that $m$ is minimal for this property.
Let $x$ be the label successor of $b_m$~in~ $\pi$.
  Then $x$ is still to the right of $b_m$ in $\psi(\pi)$, and this
  implies that its label has decreased too. By the choice of $m$, the
  letter $x$ does not belong to \B . Let $x_r \cdots x_1$ be the
  successor sequence of $x$ in $\pi$, with $x_r=x$. By Lemma~\ref{new-lemma}, none
  of the $x_i$ are in \B. 
Consequently, the order of the $x_i$ is not changed during the shift,
  so the label of~$x$ cannot have decreased, a contradiction. Thus the
  label of $b_m$ cannot decrease.
 \end{proof}

\begin{prop}[\bf{The prefix of \A }]\label{a-before-int}
  Assume \A\ and \B\ do not coincide.
  Assume $a_k, \ldots , a_{i+1}$ do not belong to \B , with $0\le i
  \le k$. The \aseq\ of $\BWX(\pi)$ begins with $a_k \cdots a_{i+1}$
  and even with $a_k \cdots a_{i}$ if $i >0$.
\end{prop}
\begin{proof}
  We first show that $a_k$ is the first letter of $\A(\BWX(\pi))$.
  Suppose not.  Let $x$ be the first letter of the new \aseq .  Then
  $x$ has label $k$ in $\BWX(\pi)$ and is smaller than $a_k$, since
  $a_k$ still has label $k$ in $\psi(\pi)$, by
  Lemma~\ref{a-labels-no-change}.  Since $x$ was already smaller than
  $a_k$ in $\pi$, it means that the label of $x$ has changed during the
  \B-shift (otherwise it would have been the starting point of the
  original \aseq).
  By Lemma~\ref{incr-label},
  the label of $x$ has actually decreased.
In other words, the label of $x$ is larger than $k$ in $\pi$.
  
  But then the successor sequence of $x$ in $\pi$ must contain a
  letter with label $k$, and this letter is smaller than $x$ and hence
  smaller than $a_k$, which contradicts the choice of $a_k$. Thus the
  first letter of the new \aseq\ is $a_k$.

  \smallskip We now prove that no letter
  can be the first (leftmost) letter that replaces one of the letters
  $a_{k-1},\ldots,a_{i}$ in the new \aseq.  Assume that the \aseq\ of
  $\psi(\pi)$ starts with $a_k\cdots a_{p+1}x$, with $i\le p<k$ and
  $x\not = a_p$. 
  Then $x$ has label $p$ in $\psi(\pi)$, and $a_p$ has label $p$
  as well (Lemma~\ref{a-labels-no-change}). Since $x$ is chosen in
  $\A(\psi(\pi))$ instead of $a_p$, this means that
$a_k, \ldots, a_{p+1}, x, a_p$ come in this order in
  $\psi(\pi)$, and that $x<a_p$. Let us prove that the letters $a_k,
  \ldots, a_{p+1}, x, a_p$    also come in this order in 
  $\pi$. Since $a_k, \ldots, a_{p+1}$ do not belong to \B, they cannot
  have moved during the shift, so it is clear that $x$ follows
  $a_{p+1}$ in $\pi$. Moreover, $x$ must precede $a_p$ in $\pi$,
  otherwise we would have $p=\ell(a_p)>\ell(x)$ in $\pi$,
  contradicting Lemma~\ref{incr-label}.

Thus $a_k, \ldots, a_{p+1}, x, a_p$ come in this order in
  $\pi$, and Lemma~\ref{incr-label} implies that
  $\ell(x)\ge p$ in $\pi$. By definition of the \aseq, $\ell(x)$
  cannot be equal to $p$. Hence $\ell(x)>p$, which forces
  $\ell(a_{p+1})>p+1$, a contradiction.

  Since no letter can be the \emm first, letter replacing one of
  $a_{k-1},\ldots,a_{i+1}, a_i$ in the new \aseq, these letters form the
  prefix of the new \aseq.
\end{proof}

The example presented at the beginning of this section shows that the
next letter of the \aseq, namely $a_{i-1}$, may not belong to the
\aseq\ after the \B-shift.

\begin{prop}[\bf{The suffix of \A}]\label{a-after-int}
  Assume \A\ and \B\ intersect but do not coincide.  Let $a_{d-1}$ be
  the first letter of $\A$ after $\A\cap\B$.  After the \B-shift, the
  \aseq\ ends with $a_{d-1} \cdots a_1$.
\end{prop}
\begin{proof}
  Observe that the existence of $a_{d-1}$ follows from
  Proposition~\ref{a-more-after}.
  
  Most of the proof will be devoted to proving that $a_{d-1}$ still
  belongs to the \aseq\ after the \B-shift.  Suppose not.  Let
  $a_m=b_p$ be the rightmost letter of $\A\cap\B$ that still belongs
  to the \aseq\ after the \B-shift (such a letter does exist, by
  Proposition~\ref{a-before-int}).  Let $a_d=b_e$ be the rightmost
  letter of $\A\cap\B$. (Note that $d-e=m-p$.) The \aseq\ of
  $\psi(\pi)$ ends with $a_m x_{m-1}\cdots x_d x_{d-1} y_{d-2} \cdots
  y_1$, with $\ell(x_j)=j$, and $x_j\not = a_j$ for $m-1 \ge j \ge
  d-1$. Let us prove that none of the $x_j$ were in the original \bseq
  .  If $x_j$ were in the original \bseq, its label in $\pi$ would
  have been $j$ (Lemma~\ref{no-b-change}). But for $m-1\ge j \ge d$,
  the only letter of $\B(\pi)$ having label $j$ is $a_j$, and by
  Lemma~\ref{b-after-int}, no letter in $\B(\pi)$ has label $d-1$.
  Thus the $x_j$ cannot have been in $\B(\pi)$. This guarantees that
  they have not moved during the \B-shift. Moreover, since they are
  smaller than $b_k$, their labels cannot have increased during the
  shift (Lemma~\ref{incr-label}).
  
  Let us prove that for $m-1\ge h \ge d-1$, 
the letter $x_h$ precedes
  $a_{h}$ in $\psi(\pi)$. We proceed by decreasing induction on $h$.
  First, $a_m$ belongs to $\A(\psi(\pi))$ by assumption, the letters
  $x_{m-1}$ and $ a_{m-1}$ are to its right and have the same label,
  and $x_{m-1}$ is chosen in the \aseq\ of $\BWX(\pi)$, which implies
  that it precedes $ a_{m-1}$.  Now assume that $x_h$ precedes $a_h$
  in $\psi(\pi)$, with $m-1\ge h \ge d$. The letter $x_h$ belongs to
  $\A(\psi(\pi))$, the letters $x_{h-1}$ and $ a_{h-1}$ are on its
  right and have the same label, and $x_{h-1}$ is chosen in the new
  \aseq, which implies that it precedes $a_{h-1}$.  Finally, $x_{d-1}$
  precedes $ a_{d-1}$ in $\psi(\pi)$, and is smaller than it.

  Let us focus on $x_{d-1}$. Assume first that it is to the right of
  $a_d$ in $\pi$. Since $\ell(x_{d-1})\ge d-1$ in $\pi$, there is a
  letter $y$ in the successor sequence of $x_{d-1}$ that has label
  $d-1$ and is smaller than $a_{d-1}$, which contradicts the choice of
  $a_{d-1}$ in the original \aseq .
  
  Thus $x_{d-1}$ is to the left of $a_d$ in $\pi$, and
  hence to the left of $b_{e-1}$.  The sequence $ b_{p+1}x_{m-1}\cdots
  x_{d-1} $ is a decreasing sequence of $\pi$ of the same length as
  $b_{p+1}b_{p}\cdots b_{e}$, and $x_{d-1}$ precedes $b_{e-1}$. By
  Lemma~\ref{key-lemma}, this implies that
  $\ell(x_{d-1})<\ell(b_{e-1})$. But $\ell(x_{d-1})\ge d-1$, 
so that $\ell(b_{e-1})\ge d = \ell(b_e)$, which is impossible.
  
  \medskip We have established that $a_{d-1}$ belongs to the \aseq\ 
  after the \B-shift. Assume now that $a_{d-1}, a_{d-2}, \ldots ,
  a_{h}$ all belong to the new \aseq , but not $a_{h-1}$, which is
  replaced by a letter $x_{h-1}$.  
This implies that  $x_{h-1}<a_{h-1}$. 
By Lemma~\ref{incr-label}, the
  label of $x_{h-1}$ was at least $h-1$ in $\pi$. Also, $x_{h-1}$ was
  to the right of $a_h$ in $\pi$. Thus in the successor sequence of
  $x_{h-1}$ in $\pi$, there was a letter $y$, at most equal to
  $x_{h-1}$, that had label $h-1$ and was smaller than $a_{h-1}$, which
  contradicts the choice of $a_{h-1}$ in the original \aseq .
\end{proof}

\smallskip
\subsection{The composition of $\bwx$ and $\BWX$}\label{section-composition}
We have seen that the beginning and the end of the \aseq\ are
preserved after the \B-shift. By symmetry, we obtain a similar result
for the \bseq\ after the \A-shift.
\begin{cor}\label{b-begin-end}
  Assume \A\ and \B\ intersect but do not coincide. Let $a_i=b_j$ be
  the leftmost element of $\A\cap\B$ and let $a_d=b_e$ be the
  rightmost element of $\A\cap\B$.  After the \A-shift, the \bseq\ 
  begins with $b_k \cdots b_{j+1}$ and ends with $a_{d-1} b_{e-1}
  \cdots b_1$.
\end{cor}
\begin{proof}
  This follows from Propositions~\ref{a-before-int}
  and~\ref{a-after-int}, together with symmetry.  Namely, since by
  Proposition~\ref{a-before-int} the first (largest) letter of the
  intersection still belongs to the \aseq\ after the \B-shift, the
  \emm place, of the last (smallest) letter of the intersection still
  belongs to the \bseq\ after the \A-shift.  After the \A-shift, the
  letter in this place is $a_{d-1}$.  The rest of the claim
  follows directly from symmetry, together with the propositions
  mentioned.
\end{proof}

\begin{figure}[htb]
\begin{center}
  \input{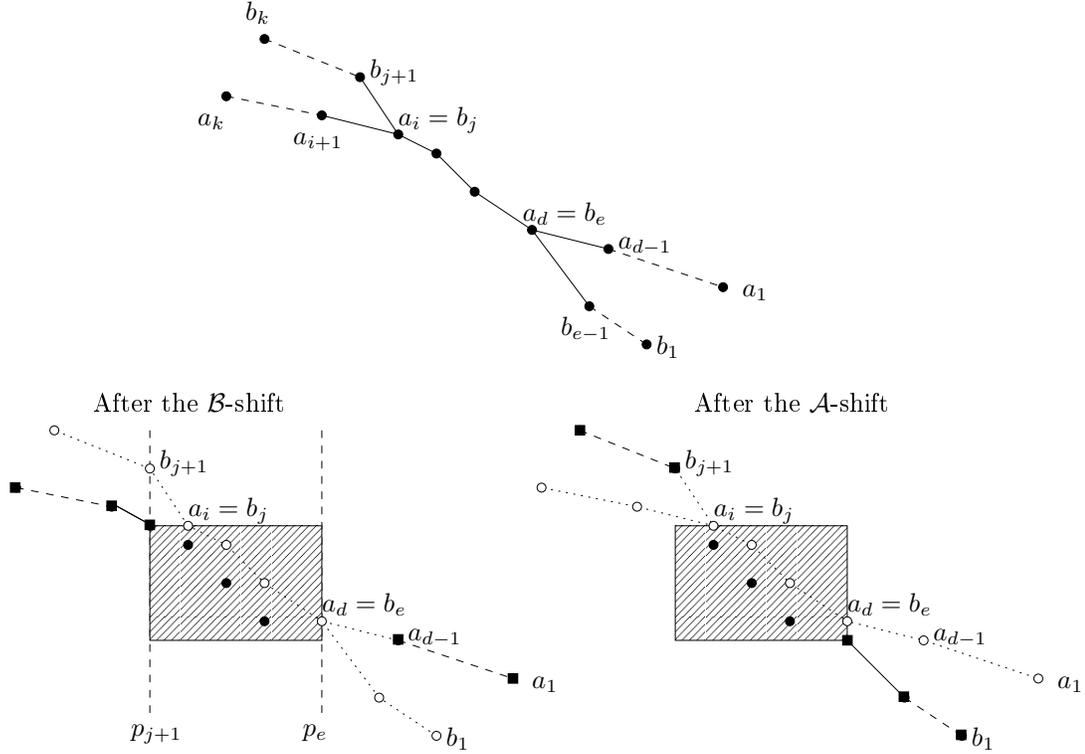}
\end{center}
\caption{The \A- and \B-sequences in $\pi$ (top), and what happens to them
  after the \B-shift (left) and the \A-shift (right). Only the black
  discs and squares belong to the permutations. The squares show
  some letters of the new \aseq\ (left) or new \bseq\ (right). The interior
  of the shaded rectangle contains the letters of~$u$.}
\label{fig-intersection}
\end{figure}

It remains to describe how the intersection of the \A- and
\B-sequences is affected by the two respective shifts.
\begin{prop}[\bf{The intersection of \A\ and \B}]\label{same-int-lemma}
  Assume \A\ and \B\ intersect but do not coincide. 
Let $a_i=b_j$ be
  the leftmost element of $\A\cap\B$ and let $a_d=b_e$ be the
  rightmost element of $\A\cap\B$. Let $a_k\cdots
 a_i x_{i-1} \cdots x_d a_{d-1}\cdots a_1$ be the \aseq\ of
  $\BWX(\pi)$. Let $b_k\cdots b_{j+1} y_{i-1} \cdots y_d
  a_{d-1}b_{e-1}\cdots b_1$ be the \bseq\ of $\bwx(\pi)$. Then
  $x_m=y_m$ for all $m$. Moreover, $x_m$ lies at the same position in
  $\BWX(\pi)$ and $\bwx(\pi)$.
\end{prop}
\begin{proof}
  First, note that the above form of the two sequences follows from
  Propositions~\ref{a-before-int},~\ref{a-after-int} and
  Corollary~\ref{b-begin-end}. Note also that if $i=d$, that is, the
  intersection is reduced to a single point, then there is nothing to
  prove.
  
  Our first objective is to prove that the sequences $\X = x_{i-1}
  \cdots x_d$ and $\Y = y_{i-1} \cdots y_d$ are the \A- and \bseq s of
  length $i-d$ of the same word (the generalization of the notion of
  \A- and \B-sequences to words with distinct letters is
  straightforward).
  
  By definition of the \aseq\ of $\BWX(\pi)$, \X\ is the smallest
  sequence of length $i-d$ (for the lexicographic order) that lies
  between $a_i$ and $a_{d-1}$ in $\BWX(\pi)$. By this, we mean that it
  lies between $a_i$ and $a_{d-1}$ both in position and in value.
  
  Let $p_m$ denote the position of $b_m$ in $\pi$. Let us show that
  \X\ actually lies between the positions $p_{j+1}$ and $p_e$
  (Figure~\ref{fig-intersection}). The first statement is clear, since
  $p_{j+1}$ is the position of $a_i$ in $\BWX(\pi)$.
  In order to prove that $x_d$ is to the left of $p_e$ in $\BWX(\pi)$,
  we proceed as at the beginning of the proof of
  Proposition~\ref{a-after-int}. We may assume $x_d\not = a_d$
  (otherwise, $x_d$ is definitely to the left of $p_e$).  Let $a_m$
  be the rightmost letter of $\A\cap\B$ that belongs to the \aseq\ 
  after the \B-shift (such a letter does exist, and $d<m\le i$).  The
  \aseq\ of $\psi(\pi)$ ends with $a_m x_{m-1}\cdots x_{d}a_{d-1}
  \cdots a_1$, with $\ell(x_j)=j$ for all $j$, and $x_j\not = a_j$ for
  $m-1 \ge j \ge d$.
  
  Let us prove, by a decreasing induction on $j\in[d,m-1]$, that the letter
  $x_j$ precedes $a_j$ for all $j$.  First, $a_m$ belongs to
  $\A(\psi(\pi))$ by assumption, the letters $x_{m-1}$ and $ a_{m-1}$
  are to its right and have the same label, and $x_{m-1}$ is chosen in
  the new \aseq, which implies that  it precedes $a_{m-1}$.
  Now assume that $x_h$ precedes $a_h$ in $\psi(\pi)$, with $m-1\ge h
  >d$.  The letter $x_h$ belongs to $\A(\psi(\pi))$, the letters
  $x_{h-1}$ and $ a_{h-1}$ are on its right and have the same label,
  and $x_{h-1}$ is chosen in the new \aseq, which implies that 
 $x_{h-1}$
  precedes $a_{h-1}$, and concludes our proof that $x_j$ precedes
  $a_j$.  
    In particular, $x_d$ is to the left of $a_d$, and hence to the left
  of the position $p_e$.
  
  We can summarize the first part of this proof by saying that \X\ is
  the smallest sequence of length $i-d$ in $\psi(\pi)$ that lies in
  position between $p_{j+1}$ and $p_e$ and in value between $a_i$ and
  $a_{d-1}$. In other words, let $u$ be the word obtained by retaining
  in $\psi(\pi)$ only the letters that lie between $p_{j+1}$ and $p_e$
  in position and between $a_i$ and $a_{d-1}$ in value. Then \X\ 
  is the \aseq\ of length $i-d$ of $u$.
  
  \smallskip By symmetry, \Y\ is the \bseq\ of length $i-d$ of the word
  $v$ obtained by retaining in $\phi(\pi)$ the letters that lie
  between $p_{j+1}$ and $p_e$ in position and between $a_i$ and
  $a_{d-1}$ in value.  But the words $u$ and $v$ actually coincide, for
  they contain
\begin{itemize}
\item[--] the letters of $\pi$ that do not belong to \A\ or \B\ and
  lie between $p_{j+1}$ and $p_e$ in position and between $a_i$ and
  $a_{d-1}$ in value. These letters keep in $\psi(\pi)$ and
  $\phi(\pi)$ the position they had in $\pi$,
\item[--] the letters 
$b_{j-1}, \ldots, b_e$, placed at positions
$p_{j}, \ldots, p_{e+1}$ (see Figure~\ref{fig-intersection}).
\end{itemize}

Observe also that $u$ does not contain any decreasing sequence of
length larger than $i-d$, because otherwise, we could use this sequence to
extend the \aseq\ of $\BWX(\pi)$. Hence we have a 
word $u$ with distinct letters,
with its \A- and \B-sequences (of length $i-d$) and we know that there
is no longer decreasing sequence in $u$. In particular, the rightmost letter
of its \bseq, $y_{i-1}$, has label $i-d$, and Lemma~\ref{bk-label}
implies that \X \ and \Y\ coincide.
\end{proof}

\begin{theorem}[{\bf Local commutation for permutations}]
\label{thm-permutations}
Let $\pi$ be a \p\  for which the \A- and \B-sequences do not
coincide.
Then $\phi(\pi)$ and $\psi(\pi)$ still contain the pattern $k\cdots 21$,
and $\bwx(\BWX(\pi))=\BWX(\bwx(\pi))$. 
\end{theorem}
\begin{proof}
The first statement follows from Lemma~\ref{a-in-order}, plus
symmetry.

  Assume first that \A\ and \B\ are disjoint. By
  Proposition~\ref{a-before-int}, the \aseq\ is unchanged after the
  \B-shift. Thus the \p\ $\phi(\psi(\pi))$ can be obtained by shifting
  \A\ and \B\ in $\pi$ in parallel. By symmetry, this is also the
  result of applying $\psi\circ \phi$ to $\pi$.
  
  Let us now assume that \A\ and \B\ intersect. Following the notation
  of Proposition~\ref{same-int-lemma}, let $a_k\cdots a_i x_{i-1}
  \cdots x_d a_{d-1}\cdots a_1$ be the \aseq\ of $\BWX(\pi)$, and let\lb
  $b_k\cdots b_{j+1} x_{i-1} \cdots x_d a_{d-1}b_{e-1}\cdots b_1$ be
  the \bseq\ of $\bwx(\pi)$.  Clearly, the only letters that 
  can move when we apply $\phi \circ\psi$ (or $\psi\circ\phi$) to $\pi$,
  are those of \A , \B\ and \X.  We need to describe at which place
  each of them ends. We denote by $p(x)$ the position of the letter
  $x$ in $\pi$ (note that $p(x)=\pi^{-1}(x)$).
  
  Let us begin with the transformation $\phi\circ \psi$. That is, the
  \B-shift is applied first. It is easy to see what happens to the
  letters that lie far away from the intersection of \A\ and \B\ 
  (Figure~\ref{fig-intersection}). During the \B-shift, the letter
  $b_k$ is sent to $p(b_1)$ and then it does not move during the
  \A-shift (it is too big to belong to the new \aseq).  Similarly,
  for $j+1 \le h <k$, and for $1 \le h <e$, the letter $b_h$ is sent
  to $p(b_{h+1})$, and then does not move. As far as the \aseq \ is
  concerned, we see that $a_h$ does not move during the \B-shift, for
  $1\le h \le d-1$ and $i+1 \le h \le k$. Then, during the \A-shift,
  $a_k$ is sent to $p(a_1)$, and the letter $a_h$ moves to
  $p(a_{h+1})$ for $1\le h <d-1$ and $i+1 \le h < k$.
  
  It remains to describe what happens to $a_{d-1}, \ldots , a_i$, and
  to the $x_h$. The letter $a_i$ moves to $p(b_{j+1})$ first, and
  then, being an element of the new \aseq , it moves to $p(a_{i+1})$.
  The letter $a_{d-1}$ only moves during the \A-shift, and it moves to
  the position of $x_{d}$ in $\psi(\pi)$. For $d\le h <i-1$, the
  letter $x_h$ moves to the position of $x_{h+1}$ in $\psi(\pi)$. The
  letter $x_{i-1}$ moves to the position of $a_i$ in $\psi(\pi)$, that
  is, to $p(b_{j+1})$. Finally, the letters $a_h$, with $d\le h <i$,
  which are not in \X\ move only during the \B-shift and end up at
  $p(a_{h+1})$.
  
  Let us put together our results: 
 When we apply $\phi\circ \psi$,
\begin{itemize}
\item[--] $x_{i-1}$ moves to $p(b_{j+1})$,
\item[--] $x_h$ moves to the position of $x_{h+1}$ in $\psi(\pi)$, for
  $d\le h <i-1$,
\item[--] $a_k$ is sent to $p(a_1)$ and $b_k$ to $p(b_1)$,
\item[--] $a_{d-1}$ moves to the position of $x_{d}$ in $\psi(\pi)$,
\item[--] the remaining $a_h$ and $b_h$ move respectively to
  $p(a_{h+1})$ and $p(b_{h+1})$.
\end{itemize}
Now a similar examination, together with the fact that each $x_h$
lies in the same position in $\psi(\pi)$ and $\phi(\pi)$
(Proposition~\ref{same-int-lemma}), shows that applying 
$\psi \circ\phi$ results exactly in the same moves.

\end{proof}

\section{Local commutation: from permutations to rook placements}
The aim of this section is to derive the local commutation for
placements (Theorem~\ref{thm-local}) from the commutation theorem for
\ps \ (Theorem~\ref{thm-permutations}).  We begin with a few simple
definitions and lemmas.

\begin{figure}[htb]
\begin{center}
  \input{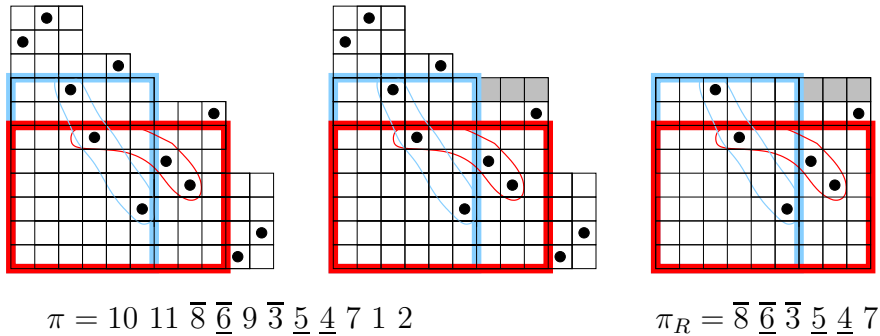}
\end{center}
\caption{Left: A placement $p$, its \A-\ and \B-sequences (for $k=3$),
  and the rectangles $R_\A$ and $R_\B$. Center: The placement $p^R$.
  Right: The placement $p_R$ and the corresponding subsequence of
  $\pi$.}
\label{fig-3placements}
\end{figure}

A \emm corner cell, $c$ of a Ferrers shape $\lambda$ is a cell such
that $\lambda \setminus \{c\}$ is still a Ferrers shape. If $p$ is a
placement on $\lambda$ containing $k\cdots 21$, with \aseq \ 
$a_k\cdots a_1$, then the \emm \A-rectangle, of $p$, denoted by
$R_\A$, is the largest rectangle of $\lambda$ whose top row contains
$a_k$. Symmetrically, the \emm \B-rectangle, of $p$, denoted by
$R_\B$, is the largest rectangle of $\lambda$ whose rightmost column
contains $b_1$ (where $b_k\cdots b_1$ is the \bseq\ of $p$). By
definition of the \A- and \B-sequences, $R_\B$ is at least as high,
and at most as wide, as $R_\A$. See the leftmost placement of
Figure~\ref{fig-3placements} for an example.

In the following lemmas, $p$ is supposed to be a placement on the
board $\lambda$, containing the pattern $k\cdots 21$.
\begin{lemma}
  Let $c$ be a corner cell of $\lambda$ that does not contain a dot
  and is not contained in $R_\A$. Let $q$ be the placement obtained by
  deleting $c$ from $p$. Then the \A-sequences of $p$ and $q$ are the
  same.
\end{lemma}
\begin{proof}
  After the deletion of $c$, the sequence $a_k\cdots a_1$ remains an
  occurrence of $k\cdots 21$ in $q$. Since the deletion of a cell
  cannot create new occurrences of this pattern, $a_k\cdots a_1$
  remains the smallest occurrence for the lexicographic order.
\end{proof}

\begin{lemma}\label{add-cells}
  Adding an empty corner cell $c$ to a row located above $R_\A$ does
  not change the \aseq . By symmetry, adding an empty corner cell 
  to a column located to the right of $R_\B$ does not change the \bseq.
\end{lemma}
\begin{proof}
  Assume the \aseq\ changes, and let $\A'=a'_k\cdots a'_1$ be the
  \aseq\ of the new placement $q$. Observe that $a_k\cdots a_1$ is
  still an occurrence of $k\cdots 21$ in~$q$.  By the previous lemma,
  $c$ belongs to $R_{\A'}$, the \A-rectangle of $q$. However, by
  assumption, $c$ is above $R_\A$.  This implies that the top row of
  $R_{\A'}$ is higher than the top row of $R_\A$, so that $a'_k$ is
  higher (that is, larger) than $a_k$. This contradicts the definition
  of the \aseq \ of $q$.
\end{proof}

\noindent
{\bf Remark.} The lemma is not true if the new cell is not added above
$R_\A$, as shown by the following example, where $k=2$. 
The \aseq \ is shown with black disks.\\
\begin{center}
  \begin{picture}(0,0)%
\includegraphics{counter-example.pstex}%
\end{picture}%
\setlength{\unitlength}{3947sp}%
\begingroup\makeatletter\ifx\SetFigFont\undefined%
\gdef\SetFigFont#1#2#3#4#5{%
  \reset@font\fontsize{#1}{#2pt}%
  \fontfamily{#3}\fontseries{#4}\fontshape{#5}%
  \selectfont}%
\fi\endgroup%
\begin{picture}(1824,474)(439,-2023)
\end{picture}%

\end{center}
\smallskip

Let $R$ be the smallest rectangle containing both $R_\A$ and $R_\B$.
It is possible that $R$ is not contained in $\lambda$. Let $p^R$ be
the placement obtained by adding the cells of $R \setminus \lambda$ to
$p$. The previous lemma implies the following corollary, illustrated
by the central placement of Figure~\ref{fig-3placements}.
\begin{cor}\label{lemma-R}
  The placements $p$ and $p ^R$ have the same \aseq \ and the same
  \bseq .
\end{cor}
\begin{proof} All the new cells are above $R_\A$ and to the right of
  $R_\B$.
\end{proof}

In what follows, the definitions of the \A- and \B-sequences, and of
the maps $\psi$ and $\phi$, are extended in a straightforward manner
to \emm partial, rook placements (some rows and columns may contain no
dot). We extend them similarly to words with distinct letters.
\begin{lemma}\label{dot-deletion}
  Let $p$ be a partial rook placement containing the pattern $k\cdots 21$. 
If we delete a row located above $a_k$,  the \aseq\ will not   change.
A symmetric statement holds for the deletion of a column located to
the right of $b_1$.
\end{lemma}
\begin{proof}
  The sequence $a_k\cdots a_1$ is still an occurrence of $k\cdots 21$
  in the new placement, and deleting a row cannot create a new
  occurrence of this pattern.
\end{proof}

\begin{prop}\label{board-perm}
Let $\pi$ be the permutation associated with a placement $p$
  containing $k\cdots 21$. There exists a
  subsequence of $\pi$ that has the same \A- and \B-sequences as~$p$.
  One such subsequence is $\pi_R$, the subsequence of $\pi$
  corresponding to the dots contained in $R$.
\end{prop}
\begin{proof}
  By Corollary~\ref{lemma-R}, we can assume that $R$ is included in
  $\lambda$. By Lemma~\ref{dot-deletion}, we can assume that
  $\lambda=R$, which concludes the proof.
\end{proof}
We shall denote by $p_R$ the (partial) placement obtained from $p^R$
by deleting all rows above $R$ and all columns to the right of $R$
(third placement in Figure~\ref{fig-3placements}).

\begin{lemma}\label{counting-lemma}
  Let $i\le j <k$. In $\psi(p)$, the maximum length of a decreasing
  sequence starting at $b_j$ and ending at $b_i$ is $j-i+1$.  One such
  sequence is of course $b_j b_{j-1} \cdots b_i$.
\end{lemma}
\begin{proof}
  Clearly, it suffices to prove the statement under the assumption
  that $b_j$ and $b_i$ are the only letters in the sequence that are
  shifted elements of the \bseq\ of $p$, which we now assume.

  Suppose that there exists in $\psi(p)$ a longer decreasing sequence, of
  the form $b_j x_j x_{j-1}\cdots x_{i+1} b_i$, where the $x$'s do not
  belong to the \bseq \ of $p$. Then $b_k \cdots b_{j+1} x_j\cdots
  x_{i+1} b_i\cdots b_1$ is an occurrence of the pattern $k\cdots 21$
  in $p$.  The fact that $x_{i+1}$ comes before $b_i$ in $\psi(p)$
  means that $x_{i+1}$ precedes $b_{i+1}$ in $p$. This contradicts the
  construction of the \bseq\ of $p$ (Lemma~\ref{dfn-psi}).
\end{proof}

The following proposition is the last technical difficulty we meet in
the proof of the commutation theorem.

\begin{prop}\label{last-lemma}
  Assume the \A- and \B-sequences of $p$ do not coincide. Then 
 $\psi(p)$ contains the pattern $k\cdots 21$, and its
  \aseq\ begins with $a_k$.
\end{prop}
\begin{proof}
 Let $R_\A$ and $R_\B$  denote the \A- and \B-rectangles of $p$. They
  form sub-boards of   $\lambda$.  Let $R$ be the smallest
  rectangle containing  $R_\A$ and $R_\B$.

Let us first prove that there exists in $\psi(p)$ an occurrence
  of $k\cdots 21$ starting with~$a_k$.  
  First, since $p$ and $p_R$ have the same \bseq\ 
  (Proposition~\ref{board-perm}), the map $\psi$ acts in the same way
  on these two placements. This means that $\psi(p_R)$ can be obtained
  from $\psi(p)$ by deleting the rows above $R_\B$ and to the right of
  $R_\A$, and by adding the cells of $R\setminus\lambda$.
  Then, by Proposition~\ref{a-before-int}, $\psi(p_R)$ contains an
  occurrence of $k\cdots 21$ starting with $a_k$, namely, the \aseq\ 
  of $\psi(p_R)$. These dots are all contained in $R_\A$, and so they
  form, in $\psi(p)$ also, an occurrence of $k\cdots 21$ starting with
  $a_k$.

  Now let $x_k \cdots x_1$ be the \aseq\ of $\BWX(p)$, and assume
  that $x_k\not = a_k$ (which implies that $x_k <a_k$).  We will
  derive from this assumption a contradiction, which will 
complete the   proof.

  If none of the values $x_j$ were in $\B(p)$, then they would form an
  occurrence of $k\cdots 21$ in $p$, which would be smaller than $a_k
  \cdots a_1$, a contradiction. Hence at least one of the $x_j$ is in
  $\B(p)$. Let $x_{\ell}=b_m$ (resp.~$x_{i+1}=b_{j}$) be the leftmost
  (resp.~rightmost) of these. Then $x_k, \ldots , x_{\ell+1}$ and
  $x_i,\ldots,x_1$ are in the same places in $p$ as in $\BWX(p)$.
  
  \medskip\noindent We consider two cases:

  \noindent
  {\bf Case 1:} Suppose first that one of the $x_r$, for $1\le r \le
  i$, lies ``above'' the \bseq\ in $p$.
  By this we mean that there exists an $s$ such that $b_s<x_r$ and
  $b_s$ precedes $x_r$ in $p$.  Let $r\le i$ be maximal such that
  $x_r$ satisfies this condition.  Let $s$ be maximal such that $b_s$
  satisfies this condition for $x_r$. Clearly, $s<k$, because $b_s<x_r
  < x_k <a_k<b_k$.
  
  The maximality of $s$ implies that $b_{s+1}>x_r$.  
In fact, $b_{s+1}$ is the smallest element of \B\ that is larger than
$x_r$.
Consider, in $p$,  the decreasing sequence
  $$
  x_k \cdots x_{\ell+1} b_m \cdots b_{s+1} x_r \cdots x_1.
  $$
  It is an occurrence of a decreasing pattern, which, given that
  $x_k<a_k$, cannot be as long as $a_k \cdots a_1$. That is,
\begin{equation}
\label{ineq}
k-\ell+m-s+r<k.
\end{equation} 

Assume for the moment that $r<i$.  By maximality of $r$, we know that
$x_{r+1}$ precedes $b_{s}$ in $p$.  Let us show that it actually
precedes $b_{s+1}$ (and thus precedes $b_s$ in $\psi(p)$).  If not,
$x_{r+1}$ lies between $b_{s+1}$ and $b_s$. But $x_{r+1}>b_s$, since
$x_{r}>b_s$, and $x_{r+1}<b_{s+1}$ by maximality of $r$. Thus
$x_{r+1}$ lies between $b_{s+1}$ and $b_s$ in position and in value,
which contradicts the definition of the \bseq\ of $p$.  Hence
$x_{r+1}$ precedes $b_{s+1}$, so the sequence $ x_\ell\cdots x_{r+1}
b_s$ in $\psi(p)$ is decreasing and has $\ell-r+1$ elements.  But
this sequence has  $x_{\ell}=b_m$ and $b_s$ as its endpoints, so, by
Lemma~\ref{counting-lemma}, it has at most $m-s+1$ points.  In other
words, $\ell-r+1\le m-s+1$, or $\ell -r \le m-s$,
contradicting~(\ref{ineq}).

Now if $r=i$, we have $s<j$ (since $b_s<x_i$ and $b_j>x_i$).  The
sequence $ x_\ell\cdots x_{i+1} $ in $\psi(p)$ is decreasing and has
$\ell-i$ elements. But this sequence has as $x_{\ell}=b_m$ and $b_j$
as its endpoints, so, by Lemma~\ref{counting-lemma}, it has at most
$m-j+1$ points.  In other words, $\ell-i\le m-j+1$.  But, since $s<j$,
this contradicts~(\ref{ineq}).
  
\medskip
  \noindent
  {\bf Case 2:} We now assume that for each $x_r$ among
  $x_i,\ldots,x_1$ there is no $s$ such that $b_s<x_r$ and $b_s$
  precedes $x_r$ in $p$.

  Lemma~\ref{counting-lemma}, applied to the subsequence $b_m=x_\ell,
  x_{\ell-1}, \ldots, x_{i+1}=b_j$ of $\psi(p)$, implies that
  $\ell-i\le m-j+1$.
That is, $i-j\ge \ell-m-1$. Now,
  $$
  b_k \cdots b_{j+1} x_i \cdots x_1
  $$
  is a decreasing sequence in $p$ of length $k-j+i\ge k+\ell-m-1$.
  At most $k-1$ of its elements can precede $b_1$, for else $b_1$
  could not be the rightmost letter of $\B(p)$.  Hence, since $b_1$
  itself does not occur in this sequence, at least $\ell-m$ of its
  elements must be preceded by $b_1$, that is, $x_{\ell-m},\ldots ,
  x_1$ all lie to the right of $b_1$.  Recall that none of the letters
  $x_i,\ldots,x_1$
  are to the right of \emm and, above any $b_s$, so $x_{\ell -m} ,
  \ldots , x_1$
  must be smaller than $b_1$.  But then
  $$
  x_k \cdots x_{\ell+1} b_m\cdots b_1 x_{\ell-m} \cdots x_1
  $$
  is an occurrence of the pattern $k\cdots 21$ in $p$, with
  $x_k<a_k$, which contradicts the definition of the \aseq .
\end{proof}

\noindent
We are finally ready for a proof of the local commutation theorem,
which we restate.

\medskip
\noindent
{\bf Theorem (same as Theorem \ref{thm-local}).}
{ \em{Let $p$ be a placement for which the \A- and \bseq s do not
  coincide. 
Then $\phi(p)$ and $\psi(p)$ still contain the pattern $k\cdots 21$, and}}
  $$
  \bwx(\BWX(p))=\BWX(\bwx(p)).
  $$

\begin{proof} 
As above, let $R$ be the smallest rectangle containing $R_\A$ and
$R_\B$. 
The first statement follows from Proposition~\ref{last-lemma} and symmetry.

We want to prove that the map
$\phi\circ\psi$ acts in the same way on the placements $p$, $p^R$ and
$p_R$. If we prove this, then, by symmetry, the same holds for the map
$\psi\circ\phi$. But the commutation theorem for \ps\ 
(Theorem~\ref{thm-permutations}) states that $\phi(\psi(p_R))=
\psi(\phi(p_R))$. Thus $\phi(\psi(p))= \psi(\phi(p))$, and we will be
done.

By Corollary~\ref{lemma-R} and Proposition~\ref{board-perm}, the
placements $p$, $p^R$ and $p_R$ have the same \bseq. Consequently,
$\psi$ acts in the same way on these three placements. In other words,
\begin{itemize}
\item[--] $\psi(p^R)$ is obtained by adding to $\psi(p)$ the cells of
  $R\setminus \lambda$; we summarize this by writing
  $\psi(p^R)=\psi(p)^R$,
\item[--] $\psi(p_R)$ is obtained by deleting from $\psi(p^R)$ the
  rows above $R$ and the columns to the right of $R$.
\end{itemize}
It only remains to prove that $\psi(p)$, $\psi(p^R) $ and $\psi(p_R)$
have the same \aseq.

By Proposition~\ref{last-lemma}, the \aseq\ of $\psi(p)$ starts with
$a_k$. This means that the \A-rectangle of $\psi(p)$ coincides with
the \A-rectangle of $p$.
Hence Lemma~\ref{add-cells}, applied to $\psi(p)$, implies that
$\psi(p)$ and $\psi(p)^R$ have the same \aseq. But
$\psi(p)^R=\psi(p^R)$, so that $\psi(p)$ and $\psi(p^R)$ have the same
\aseq.
The \aseq\ of $\psi(p^R)$, being contained in the \A-rectangle of
$p$, is contained in $R$.  By
Lemma~\ref{dot-deletion}, the \aseq s of $\psi(p_R)$ and $\psi(p^R)$
coincide.
  \end{proof}

\bigskip \bigskip

\noindent{\bf Acknowledgements.} We thank Yves M\'etivier and G\'erard Huet for
references on the confluence of rewriting systems, and Olivier Guibert
for interesting discussions on pattern avoiding involutions.

\bibliographystyle{plain} 
\bibliography{permutations}

\begin{thebibliography}{10}

\bibitem{babson-einar}
E.~Babson and E.~Steingr{\'{\i}}msson.
\newblock Generalized permutation patterns and a classification of the
  {M}ahonian statistics.
\newblock {\em S\'em. Lothar. Combin.}, 44:Art. B44b, 18 pp. (electronic),
  2000.

\bibitem{babson-west}
E.~Babson and J.~West.
\newblock The permutations {$123p\sb 4\cdots p\sb m$} and {$321p\sb 4\cdots
  p\sb m$} are {W}ilf-equivalent.
\newblock {\em Graphs Combin.}, 16(4):373--380, 2000.

\bibitem{bwx}
J.~Backelin, J.~West, and G.~Xin.
\newblock Wilf-equivalence for singleton classes.
\newblock In H.~Barcelo and V.~Welker, editors, {\em Proceeedings of the $13$th
  Conference on Formal Power Series and Algebraic Combinatorics}, pages 29--38,
  Arizona State University, May 2001.
\newblock To appear in {\em Adv. in Appl. Math.}

\bibitem{italiens}
E.~Barcucci, A.~Del~Lungo, E.~Pergola, and R.~Pinzani.
\newblock Some permutations with forbidden subsequences and their inversion
  number.
\newblock {\em Discrete Math.}, 234(1-3):1--15, 2001.

\bibitem{bona}
M.~B{\'o}na.
\newblock Exact enumeration of {$1342$}-avoiding permutations: a close link
  with labeled trees and planar maps.
\newblock {\em J. Combin. Theory Ser. A}, 80(2):257--272, 1997.

\bibitem{mbm-motifs}
M.~Bousquet-M\'elou.
\newblock Four classes of pattern-avoiding permutations under one roof:
  Generating trees with two labels.
\newblock {\em Electron. J. Combin.}, 9(2):Research Paper 19, 2003.

\bibitem{claesson-mansour}
A.~Claesson and T.~Mansour.
\newblock Enumerating permutations avoiding a pair of {B}abson-{S}teingr\'\i
  msson patterns.
\newblock {\em Ars Combinatoria}, to appear.

\bibitem{gessel-symmetric}
I.~Gessel.
\newblock Symmetric functions and {P}-recursiveness.
\newblock {\em J. Combin. Theory Ser. A}, 53(2):257--285, 1990.

\bibitem{gessel-weinstein}
I.~Gessel, J.~Weinstein, and H.~S. Wilf.
\newblock Lattice walks in {${\bf Z}\sp d$} and permutations with no long
  ascending subsequences.
\newblock {\em Electron. J. Combin.}, 5(1):Research Paper 2, 11 pp., 1998.

\bibitem{gouyou-tableaux}
D.~Gouyou-Beauchamps.
\newblock Standard {Y}oung tableaux of height {$4$} and {$5$}.
\newblock {\em European J. Combin.}, 10(1):69--82, 1989.

\bibitem{guibert-these}
O.~Guibert.
\newblock {\em Combinatoire des permutations \`a motifs exclus, en liaison avec
  mots, cartes planaires et tableaux de Young}.
\newblock PhD thesis, LaBRI, Universit\'e Bordeaux 1, 1995.

\bibitem{guibert-vexillaires}
O.~Guibert, E.~Pergola, and R.~Pinzani.
\newblock Vexillary involutions are enumerated by {M}otzkin numbers.
\newblock {\em Ann. Comb.}, 5(2):153--174, 2001.

\bibitem{huet}
G.~Huet.
\newblock Confluent reductions: abstract properties and applications to term
  rewriting systems.
\newblock {\em J. Assoc. Comput. Mach.}, 27(4):797--821, 1980.

\bibitem{jaggard}
A.~D. Jaggard.
\newblock Prefix exchanging and pattern avoidance for involutions.
\newblock {\em Electron. J. Combin.}, 9(2):Research Paper 16, 2003.

\bibitem{kitaev-mansour}
S.~Kitaev and T.~Mansour.
\newblock A survey of certain pattern problems.
\newblock Submitted, 2004.

\bibitem{krattenthaler}
C.~Krattenthaler.
\newblock Permutations with restricted patterns and {D}yck paths.
\newblock {\em Adv. in Appl. Math.}, 27(2-3):510--530, 2001.

\bibitem{regev}
A.~Regev.
\newblock Asymptotic values for degrees associated with strips of {Y}oung
  diagrams.
\newblock {\em Adv. in Math.}, 41(2):115--136, 1981.

\bibitem{schensted}
C.~Schensted.
\newblock Longest increasing and decreasing subsequences.
\newblock {\em Canad. J. Math.}, 13:179--191, 1961.

\bibitem{rodica}
R.~Simion and F.~W. Schmidt.
\newblock Restricted permutations.
\newblock {\em European J. Combin.}, 6(4):383--406, 1985.

\bibitem{stankova-4}
Z.~E. Stankova.
\newblock Forbidden subsequences.
\newblock {\em Discrete Math.}, 132(1-3):291--316, 1994.

\bibitem{stankova}
Z.~E. Stankova.
\newblock Classification of forbidden subsequences of length {$4$}.
\newblock {\em European J. Combin.}, 17(5):501--517, 1996.

\bibitem{stankova-west}
Z.~E. Stankova and J.~West.
\newblock A new class of {W}ilf-equivalent permutations.
\newblock {\em J. Algebraic Combin.}, 15(3):271--290, 2002.

\bibitem{viennot}
G.~Viennot.
\newblock Une forme g\'eom\'etrique de la correspondance de
  {R}obinson-{S}chensted.
\newblock In {\em Combinatoire et repr\'esentation du groupe sym\'etrique
  (Actes Table Ronde CNRS, Univ. Louis-Pasteur, Strasbourg, 1976)}, pages
  29--58. Lecture Notes in Math., Vol. 579. Springer, Berlin, 1977.

\bibitem{west-these}
J.~West.
\newblock {\em Permutations with forbidden subsequences, and stack-sortable
  permutations}.
\newblock PhD thesis, MIT, 1990.

\end{thebibliography}

\end{document}